\documentclass[11pt]{amsart}
\usepackage{epsfig}
\usepackage{amssymb}
\usepackage{amsfonts}
\usepackage{amsmath}
\usepackage{epsfig} 
\usepackage{float}
\usepackage{array}
\usepackage{xspace}
\usepackage{enumerate,comment,psfrag,subfigure}
\usepackage{color}
\usepackage{boldfonts}
\usepackage{mydef}
\def\rootfig{./FIGS}

\newcommand{\old}[1]{\textcolor{red}{\it }}

\usepackage{geometry}
\begin{document}

\title[Taylor-Couette dynamos]{Nonlinear dynamo in a short Taylor-Couette setup}

\author[C. Nore, J.-L. Guermond, R. Laguerre, J. L\'eorat, F. Luddens]{%
C. Nore$^{1,2,3}$,
J.-L. Guermond$^{4,\ddag}$,
R. Laguerre$^{1,5}$,
J. L\'eorat$^6$,
F. Luddens$^{1}$
}

\thanks{$^\ddag$On long leave from LIMSI (CNRS-UPR 3251), 
BP 133, 91403, Orsay, France. This author is supported in part by 
the National Science Foundation grants DMS-0510650, DMS-0713829.}

\address{$^1$Laboratoire d'Informatique pour la M\'ecanique et
les Sciences de l'Ing\'enieur, CNRS,
BP 133, 91403 Orsay cedex, France.}

\address{$^2$Universit\'e Paris Sud 11, D\'epartement de Physique, 91405
Orsay cedex, France.}

\address{$^3$Institut Universitaire de France}

\address{$^4$Department of Mathematics, Texas A\&M University
3368 TAMU, College Station, TX 77843-3368, USA.}

\address{$^5$Observatoire Royal de Belgique, Avenue Circulaire 3, B-1180 Bruxelles, Belgium.}

\address{$^6$Luth, Observatoire de Paris-Meudon, 
place Janssen, 92195-Meudon, France.}

\address{E-mail addresses: 
(CN)  {\tt nore@limsi.fr},
(FL) {\tt luddens@limsi.fr}
(RL) {\tt rlaguerr@ulb.ac.be},
(JL)  {\tt Jacques.Leorat@obspm.fr}, 
(JLG) {\tt guermond@math.tamu.edu} 
}

\keywords{Finite elements, Magnetohydrodynamics, Taylor-Couette dynamo action}

\subjclass{65N30, 76E25, 76W05}

\date{Draft version: \today}
  
%\begin{abstract}
%blabla
%\end{abstract}

\begin{abstract}
  It is numerically demonstrated by means of a magnetohydrodynamics
  code that a short Taylor-Couette setup with a body force can sustain
  dynamo action.  The magnetic threshold is comparable to what is
  usually obtained in spherical geometries. The linear dynamo is
  characterized by a rotating equatorial dipole.  The nonlinear regime
  is characterized by fluctuating kinetic and magnetic energies and a
  tilted dipole whose axial component exhibits aperiodic reversals
  during the time evolution.  These numerical evidences of dynamo
  action in a short Taylor-Couette setup may be useful for developing
  an experimental device.
\end{abstract}

\maketitle
 
%========================================================================
%========================================================================

%\begin{figure}
%\includegraphics[width=0.33\textwidth,bb= 250 120 610 650,clip=true]{H_H_phi_200_streamlines.}
%\end{figure}

\section{Introduction}
Still a century after Larmor suggested that dynamo action can be a
source of magnetic field in astrophysics, the exact mechanism by which
a fluid dynamo can be put in action in astrophysical bodies remains an
open challenge. In addition to the numerous analytical and numerical
studies that have been done since Larmor's work, it is only recently
that fluid dynamos have been produced
experimentally~\cite{Ga2000,StMu01,Monchaux07}. These experimental
dynamos have been helpful in particular to explore the nonlinear
saturation regime.  For instance, the dynamo produced in the Cadarache
experiment~\cite{Monchaux07} has an axial dipolar component and
exhibits polarity reversals that are not unlike those observed in
astronomical dynamos. The design of this experiment, however, has
peculiar features that distinguishes it from natural dynamos. The most
notable one is that the flow motion is induced by counter-rotating
impellers. This driving mechanism induces an unrealistic differential
rotation in the equatorial plane and produces a large turbulent
dissipation.  Even with a mechanical power injection close to 300 kW,
the magnetic Reynolds number of the flow of liquid sodium hardly
reaches $\Rm= 45$. Another peculiarity of this experiment is that
dynamo action has not yet been obtained by using blades made of steel.
The dynamo threshold has been reached at $\Rm=32$ by using blades made
of soft iron instead.  The objective of the present work is to
investigate an alternative driving mechanism that shares the
fundamental symmetry properties of natural dynamos, \ie axisymmetry
and equatorial symmetry (so-called $\text{SO}(2)$-$\text{Z}2$
symmetry).  The Taylor-Couette geometry is a natural candidate for
this purpose, since this configuration is already known to produce
dynamos both in axially periodic geometries~\cite{WB02} and in finite
vessels of large aspect ratio~\cite{GLLN09}. We examine in the present
paper the dynamo capabilities of Taylor-Couette flows in vessels of
small aspect ratio, and we compare the results obtained in this
setting with those from more popular spherical
dynamos~\cite{Dudley89}.

The paper is organized as follows.  The formulation of the problem and
the physical setting of the Taylor-Couette configuration under
consideration are described in \S\ref{sec:formulation_problem}.  The
formulation of the problem and the physical setting of the
Taylor-Couette configuration under consideration are described in
\S\ref{sec:formulation_problem}. Three types of flows are considered
in the paper and are discussed in \S\ref{sec:hydro_shortTC}.  These
flows are the standard Taylor-Couette flow driven by viscous stresses,
a manufactured Taylor-Couette flow, and an optimized flow driven by a
body force that models rotating blades attached to the lids.  Two
kinematic dynamo configurations are investigated in
\S\ref{sec:dynamo}. It is found that the poloidal to toroidal ratio of
the velocity field generated by viscous driving only (standard
Taylor-Couette) is not large enough to generate a dynamo at $\Rm \le 200$.  Dynamo
action is obtained by using the strengthen Taylor-Couette flow and the
forced Taylor-Couette flow.  In both cases the poloidal to toroidal
ratio of the velocity field is close to one.  A nonlinear dynamo
obtained with the forced Taylor-Couette flow is described in
\S\ref{sec:nl_dynamo}. In the early linear phase of the dynamo, the
magnetic field at large distance is dominated by an equatorial
rotating dipole.  In the established nonlinear regime, an axial
axisymmetric component of the magnetic dipole is excited and exhibits
aperiodic reversals.  Concluding remarks are reported in
\S\ref{Sec:remarks}.

\section{Formulation of the problem}
\label{sec:formulation_problem}
\subsection{The physical setting}
\label{sec:model}
We consider an incompressible conducting fluid of constant
density $\rho$ and constant kinematic viscosity $\nu$. 
This fluid is contained between two coaxial cylinders of height $L_z$.
The radius of the inner cylinder is $R_i$ and that of the outer one is
$R_o$.  The inner cylinder is composed of a solid conducting material.
The inner cylindrical wall and the top and bottom lids corotate at
angular velocity $\Omega_i$.  The outer cylindrical wall is
motionless. The inner solid core may rotate or not, \ie the inner core
and the inner cylindrical wall may have different angular
velocities. The conductivity of the fluid and inner solid is assumed
to be constant and is denoted $\sigma_0$.  The magnetic permeability
$\mu_0$ is assumed to be constant in the entire space.

Let $\calU$ be a reference velocity scale yet to be defined.  We then
consider the following reference scales for length, $\calL = R_o - R_i$,
magnetic field, $\calH =\calU\sqrt{\rho/\mu_0}$, and pressure,
$\calP=\rho \calU^2$. The non-dimensional parameters of the system are
the kinetic Reynolds number, $\Re$, the magnetic Reynolds number,
$\Rm$, the radius ratio $\eta$, and the aspect ratio $\Gamma$:
\begin{equation}
\Re= \frac{\calU \calL}{\nu}, 
\quad 
\Rm = \mu_0 \sigma_0 \calU \calL,
\quad
\eta = \frac{R_i}{R_o},
\quad
\Gamma = \frac{L_z}{\calL}.
\end{equation}

% This then leads us to consider an effective kinetic Reynolds number,
% $\Rerms$, and an effective magnetic Reynolds number, $\Rmrms$, as
% follows: \Question{JLG: We have not used this notation}
% \begin{equation}
% \Rerms= \frac{\calU^\text{rms} \calL}{\nu}, 
% \quad 
% \Rmrms = \mu_0 \sigma_0 \calU^\text{rms} \calL.
% \end{equation}

To limit the number of geometrical parameters, we restrict ourselves
in this paper to $\eta=0.5$ and $\Gamma=2$. Abusing the notation, this
immediately implies that $R_i=1$ and $R_o=2$ in non-dimensional
units. We did not explore other aspect ratios (see for
example~\cite{Mullin91,Abshagen04,Lopez06,Mullin08} for short aspect ratios and different
angular velocities). 
% To the best of your knowledge, configurations with inner and outer
% walls both rotating are rarely studied~\cite{Mullin91}.  
The conducting domain $\Omega_{c}$ is partitioned into its fluid part
enclosed between the two walls, $\Omega_{cf}$, and its solid
part enclosed in the inner cylinder, $\Omega_{cs}$. Using
non-dimensional cylindrical coordinates $(r, \theta, z)$, we have
$\Omega_{cf}= [1,2]{\times}[0,2\pi){\times}[-1,1]$ and $\Omega_{cs} =
[0,1]{\times} [0,2\pi){\times}[-1,1]$.  The conducting material is
embedded in a non-conducting region denoted $\Omega_v$, which we refer
to as the vacuum region.

The non-dimensional set of equations that we consider is written as
follows in the conducting material:
\begin{eqnarray}
  \partial_t\bu + (\bu\ADV)\bu + \GRAD p
  & = & \frac{1}{\Re}{\boldsymbol\LAP}\bu \label{eq:nsp}+ 
  (\ROT\bHc ) \CROSS \bHc + \bef_I \\ \DIV \bu & = & 0
  \label{eq:divp}\\ \partial_t\bHc -\ROT(\bu \times \bHc) & = &
  \frac{1}{\Rm}{\boldsymbol\LAP} \bHc \label{eq:ind} \\ \DIV\bHc & = & 0,
\label{eq:divh}
\end{eqnarray}
where $\bu$, $p$, and $\bHc$ are the velocity field, pressure, and
magnetic field, respectively. The magnetic field in $\Omegav$ is
assumed to derive from a harmonic scalar potential: $\bHv=\GRAD\phi$,
$\LAP\phi=0$.  The transmission conditions across the interface
separating the conducting and nonconducting material are such that the
tangent components of the magnetic and electric fields are continuous
(see \cite{GLLN05}).  

We consider three different settings: (i) The incompressible
Navier-Stokes setting ($\bHc=0$); (ii) The Maxwell or kinematic dynamo
setting; (iii) The nonlinear magnetohydrodynamics setting (MHD).  In
the Navier-Stokes setting, $\bHc$ is set to zero in the Lorentz force
and the induction equation is not solved. The source term $\bef_I$ is
an ad hoc body force that models blades fixed at the endwalls, see
\S\ref{sec:hydro_forced_TC}. When $\bef_I=0$, the viscous stress
induced by the rotating walls is the only source of momentum, see
\S\ref{sec:hydro_viscous_TC}. In the Maxwell setting, only the
induction equation is solved assuming that some ad hoc velocity field
$\bu$ is given. In the MHD setting, the full set of equations is
solved.

Since the definition of the reference velocity in similar dynamo
configurations may be different (velocity at a given point, maximal
speed in the flow, \etc), we introduce the root mean square (rms)
velocity to facilitate comparisons:
\begin{equation}
\calU^{*2} = \frac{1}{\text{vol}(\Omega_{cf})} \int_{\Omega_{cf}} \|\bu(\bx,t)\|^2 \dif\bx,
\end{equation} 
where the dimensionless fluid volume is $\text{vol}(\Omega_{cf}) =
6\pi$ in the present case.
 
\subsection{Numerical details}
The code (SFEMaNS) that we have developed solves the coupled
Navier-Stokes and Maxwell equations in the MHD limit in heterogeneous
axisymmetric domains composed of conducting and nonconducting regions
by using a mixed Fourier/Lagrange finite element technique.
Continuous Lagrange Finite elements are used in the meridian plane and
Fourier modes are used in the azimuthal direction. Parallelization is
done with respect to the Fourier modes. Continuity conditions across
interfaces are enforced using an interior penalty technique
\cite{GLLN05,GLLN09}. SFEMaNS can account for discontinuous electrical
conductivity and magnetic permeability
distributions~\cite{GAFD_Giesecke_2010b,GLLNR11}. The magnetic field in the nonconducting
regions is assumed to derive from a scalar magnetic potential, \ie the
configurations that we model are such that there is some mechanism
that ensures that the circulation of the magnetic field along any path
in the insulating medium is zero (this happens for instance when the
vacuum is simply connected).  The velocity field in $\Omega_{cf}$ and
the magnetic field in $\Omega_c$ are approximated using continuous
$\polP_2$ polynomials, and the pressure field in $\Omega_{cf}$ is
approximated using continuous $\polP_1$ polynomials.  In the vacuum
$\Omega_v$, the magnetic potential $\phi$ is approximated using
continuous $\polP_2$ polynomials. Typical characteristics of the
meshes in the meridian section of all the cases studied in this paper
are summarized in Table~\ref{tab:runs}.

\begin{table}[h!]
\begin{tabular}{|c|c|c|c|c|c|c|c|}
  \hline 
  Run & $\Delta x$  & $\Delta t$ & $np(P)$& $np(V)$ & $np(H)$ & $np(\phi)$ & M \\ \hline
  \S\ref{sec:hydro_viscous_TC} & $1/100$ & 0.025 & 5911     &23341 & - & - & 8 \\
   \S\ref{sec:hydro_forced_TC}, \S\ref{Sec:def_of_Vepsilon} & $1/100$ & 0.025 &5911     &23341 & - & - & 12 \\ 
  \S\ref{sec:kin_pol_tor_TC}, \S\ref{sec:kin_forced_TC} & $1/100$ & 0.005 & 5911     &23341 & 29821 & 14041 & 4 \\
  \S\ref{sec:nl_dynamo}& $1/100$ & 0.005 & 5911     &23341 & 29821 & 14041 & 32 \\ \hline
\end{tabular}
\vspace{\baselineskip} \\
\caption{Characteristics of the runs: 
  $\Delta x$ is the quasi-uniform meshsize in $\Omega_c$; $\Delta t$ is the timestep;
  $np(P)$ is the number of $\polP_1$ nodes for the pressure field in $\Omega_{cf}$;
  $np(V)$ is the number of $\polP_2$ nodes for the
  velocity field in $\Omega_{cf}$;
  $np(H)$ is the number of $\polP_2$ nodes for the magnetic field in $\Omega_c$;
  $np(\phi)$ is the number of $\polP_2$ nodes for the magnetic potential in $\Omega_v$.
  The numbers $np(P)$, $np(V)$, $np(H)$ refer only to the meridian section.
  The total number of grid points for each unknown $Y$ is obtained by multiplying $np(Y)$
  by 2 times the number of Fourier modes, $M$, minus one.
}
\label{tab:runs}
\end{table}
%CYL50_G2_SYM_Z_MHD.FEM
%==> Mesh_1_2_FE_2 <==
% np_lect 43369 nw_lect 6 nws_lect 3 me_lect 21630 mes_lect 454
% np  29821 me  14774 mes  272 nps  544
%
%==> Mesh_1_FE_1 <==
% np_lect 10870 nw_lect 3 nws_lect 2 me_lect 21630 mes_lect 454
% np  5911 me  11520 mes  300 nps  300
%
%==> Mesh_1_FE_2 <==
% np_lect 43369 nw_lect 6 nws_lect 3 me_lect 21630 mes_lect 454
% np  23341 me  11520 mes  300 nps  600
%
%==> Mesh_3_FE_2 <==
% np_lect 43369 nw_lect 6 nws_lect 3 me_lect 21630 mes_lect 454
% np  14041 me  6856 mes  328 nps  656

%CYL100_G2_SYM_Z.FEM
%==> Mesh_1_FE_1 <==
% np_lect 26093 nw_lect 3 nws_lect 2 me_lect 51368 mes_lect 1016
% np  26093 me  51368 mes  1016 nps  1014
%
%==> Mesh_1_FE_2 <==
% np_lect 103553 nw_lect 6 nws_lect 3 me_lect 51368 mes_lect 1016
% np  103553 me  51368 mes  1016 nps  2030

The performance of SFEMaNS has been validated on various kinematic and
nonlinear dynamo configurations. In particular, a study of two
Taylor-Couette setups using SFEMaNS is reported in~\cite{GLLN09}. In
the first case $\Gamma=4$, $\eta=0.5$, and $z$-periodicity is assumed;
in the second case $\Gamma=2\pi$, $\eta=0.5$ and the vessel is finite,
\ie no $z$-periodicity is assumed and the vessel is closed at both
ends. In both cases the inner wall rotates, but the outer wall and the
two endwalls (when present) are motionless. The self-consistent
saturated dynamo found in~\cite{WB02} in the $z$-periodic case has
been reproduced in~\cite{GLLN09}, and a new nonlinear dynamo has been
found in the finite vessel at $\Re=120, \, \Rm=240$.  The behaviors of
the $z$-periodic and finite-vessel dynamos, as observed in
\cite{GLLN09}, significantly differ.  After some transient, the
kinetic and magnetic energies of the $z$-periodic dynamo converge to a
stationary value. The final nonlinear MHD state is a steady rotating
wave resulting from the balance between the driving effect of the
viscous shear and the braking effect of the Lorentz force.  The
nonlinear dynamo action found in the finite vessel shows a different
behavior in which the spatial symmetry about the equatorial plane (or
mid-plane) of the velocity and magnetic fields plays a key role. The
dynamo is cyclic in time and the fields rotate rigidly with modulated
amplitude.  In these two cases (periodic and finite extension), the
wavelength of the magnetic eigenvector is about twice that of the
flow; as a result, the velocity field in the median plane of a
single magnetic structure is directed inwards.  This feature is shared
by the spherical kinematic dynamos studied in~\cite{Dudley89}. It is
reported in \cite{Dudley89} that the lowest critical magnetic Reynolds
number is obtained when the velocity field forms two poloidal
cells that flow inwards in the equatorial plane.
Note in passing that the two numerical experiments
reported in~\cite{GLLN09} clearly confirm that assuming periodicity or
enforcing finite boundary conditions give rise to dynamos with
fundamentally different behaviors, \ie assuming periodicity or ad hoc
boundary conditions for the sake of numerical convenience may have
nontrivial consequences. 
The series of observations above have led us to investigate more
thoroughly the Taylor-Couette configuration with aspect ratio
$\Gamma=2$.

\section{Hydrodynamic forcing}
\label{sec:hydro_shortTC}
Since a number of dynamo studies have shown that the ratio of poloidal
to toroidal speed should be close to unity to obtain the lowest
critical magnetic Reynolds number, it is important to control this
ratio.  We describe in this section the mechanisms that we use to
optimize the velocity field for dynamo action.

\subsection{Taylor-Couette flow (viscous driving only)}
\label{sec:hydro_viscous_TC}
When the aspect ratio is about 2 and the kinematic Reynolds number is
moderate, two counter-rotating poloidal cells form with a toroidal
angular velocity oriented in the same direction as that of the inner
cylinder. In order to enforce the equatorial jet to flow inwards, we
let the lids of the vessel rotate with the angular velocity of the
inner cylinder and we keep the outer cylinder motionless.  Note that
it is important to have the lids and the inner cylindrical wall of the
vessel to corotate; this makes the equatorial jet flow inwards and
makes the overall velocity filed similar to the spherical flows that
are known to yield dynamo action~\cite{Dudley89}.

We define the velocity reference scale to be
\begin{equation}
  \calU = \Omega_i R_i,
\end{equation} 
when the only source of momentum is the viscous stress at the
boundary.

\begin{figure}[h]
\includegraphics[width=0.5\textwidth]{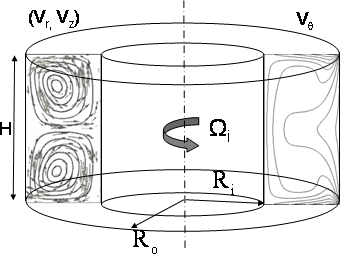}
\caption{Taylor-Couette flow $\calV_0$, $\Gamma=2$, $\Re=120$.  The
  angular velocity of the lids and the inner cylinder is $\Omega_i=1$;
  the outer cylinder is motionless.  A radial jet flows inward at the
  equator.  Represented are the poloidal flow (vectors and
  streamlines) $-0.2 \le V_r \le 0.3$, and $-0.25 \le V_z \le 0.25$,
  and the toroidal (or azimuthal) flow $0.25 \le V_\theta \le 2$
  (every 0.25).}
\label{fig:V_tch}
\end{figure}

At $\Re=120$ in the Navier-Stokes regime, the flow is steady, and
forms the expected two toroidal cells invariant under the
$\text{SO}(2)$-$\text{Z}2$ symmetry, \ie axisymmetric and symmetric
with respect to the equatorial plane, see Figure~\ref{fig:V_tch}. This
flow, henceforth generically referred to as $\calV_0$, is
characterized by its rms velocity, $V_0^*$, defined as follows:
\[
V_0^{*2}=V_{0p}^{*2}+V_{0t}^{*2}=\frac{1}{\text{vol}(\Omega_{cf})}\int_{\Omegacf}
(V_r^{2} + V_z^{2})\dif \bx +
\frac{1}{\text{vol}(\Omega_{cf})}\int_{\Omegacf} V_\theta^{2} \dif
\bx,
\]
where $V_{0p}^*$ and $V_{0t}^*$ are the rms poloidal and toroidal
velocities of the reference hydrodynamic flow, respectively, and
$\text{vol}(\Omega_{cf})=6\pi$ is the volume of the vessel. Our
computations give $V_0^{*}=0.272$; this value is significantly lower
than the maximum speed at the rim of the endwalls which is equal to
$2$. The poloidal to toroidal ratio is
$\Lambda_0=V_{0p}^{*}/V_{0t}^{*}=0.235$.

We have verified that the flow $\calV_0$ is stable with respect to
non-axisymmetric perturbations supported on the Fourier modes $m=1,
\cdots, 7$ at $\Re=120$. The velocity field $\bV_0$ together with a
sketch of the setup is shown in Figure~\ref{fig:V_tch}. This reference
flow is used in \S\ref{sec:kin_pol_tor_TC} to perform kinematic dynamo
simulations .

%\begin{figure}
%\epsfig{file=\rootfig/Ur_Re120.ps,height=0.35\textheight,bb=
%83 31 464 690,clip}\hfil
%\epsfig{file=\rootfig/Ut_Re120.ps,height=0.35\textheight,bb=
%83 31 464 690,clip}\hfil
%\epsfig{file=\rootfig/Uz_Re120.ps,height=0.35\textheight,bb=
%83 31 464 690,clip} \\
%(a) \hfil\hfil (b) \hfil \hfil(c)
%\caption{Taylor-Couette flow $\bV_{\text{VF}}^*$ with $\Gamma=2$ and $\Re=120$.
%The lids and the inner cylinder are rotating at $\Omega_i=1$ while the outer cylinder remains static.
%The radial jet is inwards at the equator.
%Represented are the radial,  azimuthal and axial velocity field components:
%$-0.2 \le V_r \le 0.3$ (every 0.05 except 0), $0.25 \le V_\theta \le 2$ (every 0.25)
%and $-0.25 \le V_z \le 0.25$ (every 0.05 except 0).}
%\label{fig:V_tch}
%\end{figure} 

\subsection{A modified Taylor-Couette flow}\label{Sec:def_of_Vepsilon}
In order to perform kinematic dynamo simulations with a velocity field
that has a poloidal to toroidal ratio that can be controlled easily,
we construct an ad hoc field based on $\calV_0$. We use the poloidal
and toroidal components of the vector field $\bV_0$ to define a
kinematic field, $\bV_\epsilon$, with a pre-assigned poloidal to
toroidal ratio as follows:
\begin{equation}
\bV_\epsilon = \frac{\epsilon}{\alpha(\epsilon)} \bV_{0p} + \frac{1}{\alpha(\epsilon)}\bV_{0t}.
\end{equation}
The normalization is done so that the rms of $\bV_\epsilon$ is the
same as that of $\bV_0$. This gives
\begin{equation}
  \alpha^2(\epsilon) = \frac{1+\epsilon^2 \Lambda_0^{2}}{1+\Lambda_0^{2}}, \qquad 
  \Lambda(\epsilon) = \epsilon \Lambda_0.
\end{equation}
Since the toroidal component of the velocity at the inner cylinder is
equal to $1/\alpha(\epsilon)$, the angular velocity of the inner wall
is $\Omega_i=1/\alpha(\epsilon)$, and this also means that the
reference velocity scale is
\begin{equation}
  \calU = \alpha(\epsilon) \Omega_i R_i. \label{def_of_U_modified}
\end{equation}

Although the vector field $\bV_\epsilon$ is not a solution of the
Navier-Stokes equations, it is nevertheless solenoidal. This flow is
henceforth generically referred to as $\calV_\epsilon$. Computations have been
done (see \S\ref{sec:kin_pol_tor_TC}) for the values of $\epsilon$
reported in Table~\ref{tab:runs_data_hydro_epsilon}. The quantity
denoted $V_{\max}$ in Table~\ref{tab:runs_data_hydro_epsilon} is the
maximum of the velocity modulus; $V_{\max}$ depends on $\epsilon$.
\begin{table}[h!]
\begin{tabular}{|c|c|c|c|c|c|c|c|c|c|c|}
\hline 
$\epsilon$         & 1     & 3    & 4    & 5    & 6    & 6.5  & 8    & 10   & 12   & 16   \\ \hline 
$\alpha(\epsilon)$ & 1     & 1.19 & 1.34 & 1.50 & 1.69 & 1.78 & 2.08 & 2.49 & 2.92 & 3.80 \\ \hline 
$\Lambda(\epsilon)$& 0.235 & 0.71 & 0.94 & 1.18 & 1.41 & 1.53 & 1.89 & 2.36 & 2.83 & 3.77 \\ \hline
$V_{\max}$ & 2.00  & 1.67 & 1.49 & 1.32 & 1.20 & 1.21 & 1.23 & 1.25 & 1.26 & 1.27 \\ \hline 
\end{tabular} 
\vspace{\baselineskip}\\
\caption{Modified Taylor-Couette flow: normalization factor $\alpha(\epsilon)$, 
poloidal to toroidal ratio $\Lambda(\epsilon)$ and maximum of the velocity modulus $V_{\max}$.}
\label{tab:runs_data_hydro_epsilon}
\end{table}

\subsection{Forced Taylor-Couette flow (viscous driving plus body
  force)}
\label{sec:hydro_forced_TC}
A number of dynamo studies have shown that the ratio of poloidal to
toroidal speed should be close to unity to obtain a low critical
magnetic Reynolds number. Viscous driving by the rotating walls yields
a value for this ratio that is not close to unity ($\Lambda_0=0.235$
at $\Re=120$, see section above). At low Reynolds numbers, the flow is
steady and axisymmetric.  It is relatively
easy to vary the relative amplitude of the toroidal component in
experimental setups by using blades fixed to the corotating endwalls
to act as centrifugal pumps. This configuration, however, is difficult
to implement in a computer code. In order to better control the
poloidal to toroidal ratio in our simulations, we have chosen to model
the toroidal driving by a body force.
The action of blades on the top and bottom lids is modeled by an ad
hoc axisymmetric divergence-less force given in dimensional form as
follows:
\begin{equation}
  \bef_I(r,z) = \begin{cases} \displaystyle \rho \frac{A}{r} \calU^2 \be_r & \text{if $0.8 \le |z| \le
      1$ and $1.2 \le r \le 1.8$} \\ 0 & \text{otherwise},  \end{cases} \label{def_calU_calV_I}
\end{equation} 
where the non-dimensional parameter $A$ has been tuned to optimize the
poloidal to toroidal ratio. Note that \eqref{def_calU_calV_I} defines
the reference velocity $\calU$.  The resulting velocity field is
denoted $\bV_I$ and the flow is generically called $\calV_I$.

\begin{figure}[h]
\subfigure[Poloidal flow, $\theta=\pi$\hspace{-0.74cm}]{
  \includegraphics[height=0.30\textwidth]{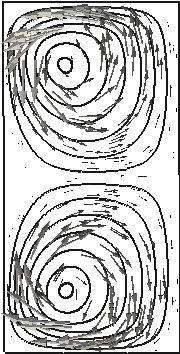}} \hfil
\subfigure[Toroidal flow, $\theta=0$\hspace{-0.78cm}]{
\includegraphics[height=0.30\textwidth]{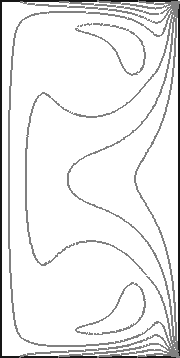}}
  \caption{Forced Taylor-Couette flow $\calV_{\text{I}}$, $\Gamma=2$,
    $\Re=120$, and $A=2.5$.  The lids and the inner cylinder rotate
    with angular speed $\Omega_i=0.55$; the outer cylinder is
    motionless. A radial jet flows inwards at the equator.
    Represented are the poloidal flow (vectors and streamlines), $-0.7
    \le V_r \le 1.4$, $-1.1 \le V_z \le 1.1$, and the toroidal (or
    azimuthal) flow, $ 0 \le V_\theta \le 1.6$ (10 levels)
% with $1 \le r \le 2, \, -1 \le z \le 1$ (the $z-$axis is on the left)
  }
  \label{fig:V_tch_forcage}
\end{figure}

We have found that using $A=2.5$ at $\Re=120$ gives
$\Lambda_I=V_{Ip}^{*}/V_{It}^{*}=1.04$ and the rms velocity is
$V_I^*=0.219$.  We have observed that the azimuthal velocity in the
vicinity of the inner radius is close to $0.55$; hence, to reduce the
viscous boundary layer at the inner wall and endwalls, we have set the
dimensionless angular velocity to $\Omega_i=0.55$. 
% Other characteristics values are summarized in
% table~\ref{tab:runs_data_hydro}.
The steady axisymmetric flow $\calV_{\text{I}}$ is shown in
Figure~\ref{fig:V_tch_forcage}, (see Figure~\ref{fig:V_tch} for a
comparison with the pure Taylor-Couette flow).
%\begin{table}[h!]
%\begin{tabular}{|c|c|c|c|}
%  \hline
%  Run & $V_{max}$ & $\Lambda$ &  $V^*$\\ \hline
%  \S\ref{sec:hydro_viscous_TC}, viscous flow $\calV_0$ & $2.00$ & 0.24 & 1.18\\
%   \S\ref{sec:hydro_forced_TC}, forced flow $\calV_{\text{I}}$ & $1.09$ & 1.04 & 0.95     \\ 
% \hline
%\end{tabular}
%\caption{Characteristics of the hydrodynamic runs: 
%  $V_{max}$ is the maximum of the velocity modulus in the fluid domain; 
%  $\Lambda=V_{p}/V_{t}$ is the poloidal to toroidal ratio;
%  $V^*$ is the r.m.s. speed.
%}
%\label{tab:runs_data_hydro}
%\end{table}
%
We have verified, by performing nonlinear Navier-Stokes simulations,
that the flow $\calV_I$, at $Re=120$, is stable with respect to
three-dimensional perturbations supported on Fourier modes up to
$m=11$. The first hydrodynamic non-axisymmetric instability occurs on
the Fourier mode $m=3$ at $\Re= 168$. The steady and axisymmetric
forced Taylor-Couette flow $\calV_{\text{I}}$ is used in
\S\ref{sec:kin_forced_TC} to perform kinematic dynamo simulations.

\subsection{Summary}
To compare the flows $\calV_0$, $\calV_I$, and $\calV_\epsilon$, we
show in Table~\ref{tab:runs_data_hydro} the following characteristics
of these three flows: rms velocity, $V^*$; maximum of the velocity
modulus in the fluid domain, $V_{\max}$; poloidal to toroidal ratio,
$\Lambda$.
\begin{table}[h!]
  \begin{tabular}{|c|c|c|c|}
    \hline
    Run & $V^*$ & $V_{\max}$ & $\Lambda$ \\ \hline
    Viscous flow $\calV_0$, \S\ref{sec:hydro_viscous_TC} 
    &  0.272 & 2.00 & 0.235 \\
    Modified flow $\calV_\epsilon$, \S\ref{Sec:def_of_Vepsilon},\S\ref{sec:kin_pol_tor_TC}
    &  0.272 & (1.20,2.00) & $0.235{\times} \epsilon$\\
    Forced flow $\calV_{\text{I}}$, \S\ref{sec:hydro_forced_TC},\S\ref{sec:kin_forced_TC},\S\ref{sec:nl_dynamo}
    & 0.219 & 1.09 & 1.04 \\
    \hline
\end{tabular}
\vspace{\baselineskip}\\
\caption{Characteristics of the flows $\calV_0$, $\calV_\epsilon$ and $\calV_{\text{I}}$:
  $V^*$ is the r.m.s. speed; $V_{\max}$ is the maximum 
  of the velocity modulus in the fluid domain; 
  $\Lambda$ is the poloidal to toroidal ratio.}
\label{tab:runs_data_hydro}
\end{table}

\section{Kinematic Dynamos}
\label{sec:dynamo}
We evaluate in this section the properties of the kinematic dynamos
generated by the flows $\calV_\epsilon$ (viscous driving) and
$\calV_I$ (viscous driving plus body force).
%JLG Commented
%\subsection{Linear regime}
%\subsubsection{Influence of the poloidal to toroidal ratio in
%  kinematic dynamo computations using the viscous forcing}}

\subsection{Parametric study of the poloidal to toroidal ratio using $\calV_\epsilon$}
\label{sec:kin_pol_tor_TC}
We investigate the dynamo properties of the manufactured flow
$\calV_\epsilon$ in the kinematic regime, see
\S\ref{Sec:def_of_Vepsilon}. The reference velocity scale is defined
in \eqref{def_of_U_modified}.  To ensure that the velocity is
continuous across the solid/fluid interface, the angular velocity of
the inner core is set to be $1/\alpha(\epsilon)$.  The conductivities
of the solid inner core and the fluid are identical.

\begin{figure}[h]
\includegraphics[width=0.5\textwidth]{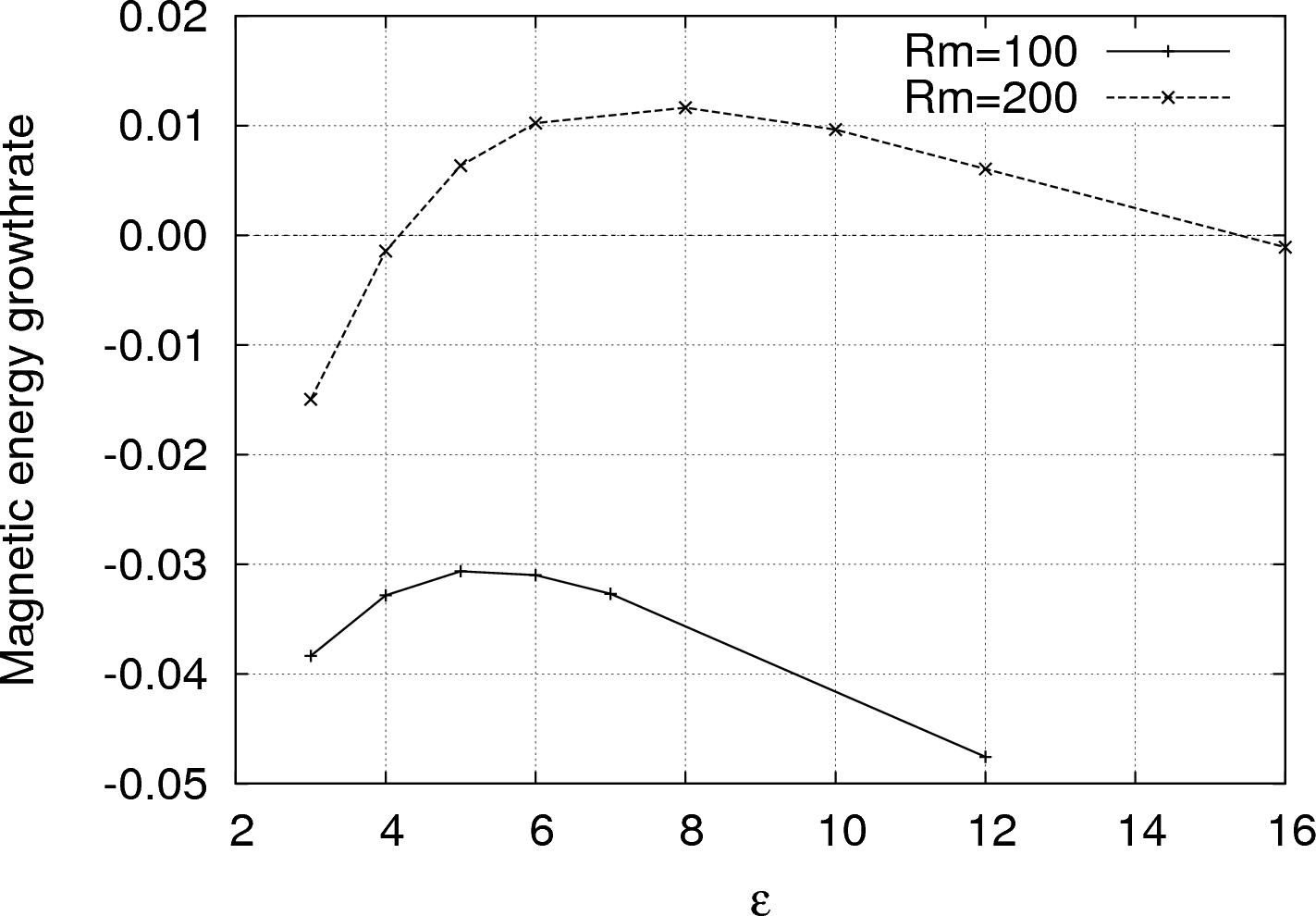}
\caption{Kinematic dynamo. Growthrate of Fourier mode $m=1$ for the
  modified Taylor-Couette flow, $\calV_\epsilon$, as a function of
  $\epsilon$ for $\Rm=100$ and $\Rm=200$; $\Gamma=2$ and
  $\Re=120$.}
\label{fig:eps_taux_Re120}
\end{figure}

%The relevant data are summarized in table~\ref{tab:runs_data_hydro}.
%
%
We perform two studies at $\Rm=100$ and $\Rm=200$ to determine the
optimal weight $\epsilon$ that gives the largest growthrate of the
dynamo action in the kinematic regime. The computations are done with
SFEMaNS in Maxwell mode. The magnetic field is initialized to some
small random values and the growth rate (\ie the real part of the
leading eigenvalue) is computed by running short time simulations for
various ratios $\epsilon \in [3,16]$ shown in
Table~\ref{tab:runs_data_hydro_epsilon}. As the vector field
$\bV_\epsilon$ is axisymmetric, the term $\ROT(\bV_\epsilon {\times}
\bHc)$ cannot transfer energy between the Fourier modes of $\bHc$, \ie
the Fourier modes are uncoupled.  The first bifurcation is of Hopf
type and the most unstable eigenvector is the Fourier mode $m=1$.  The
growthrate of the magnetic energy is reported in
Figure~\ref{fig:eps_taux_Re120}. There is no dynamo action at
$\Rm=100$. Dynamo action occurs at $\Rm=200$ in the range $4.2 <
\epsilon <15.4$, which corresponds to $ 1.0 < \Lambda(\epsilon) <
3.8$.  Note that the purely viscous driving, which corresponds
to $\epsilon=1$ and $\Lambda_0=0.235$, cannot sustain a dynamo at
$\Rm=100$ and $\Rm=200$.

%/workdir/coro/nore/TC_VARGAS/TCM/Re120_Rm200/AXE_COND/eps_8
%/workgpfs/rech/nor/rnor522/TC/AXE_COND/MXX/Re120_Rm200/Re120_Rm200_eps8/tfin_240
\begin{figure}[h]%
\begin{minipage}{0.70\textwidth}%
\subfigure[$H_r$ at $\theta=0$]{\includegraphics[width=0.3\textwidth]{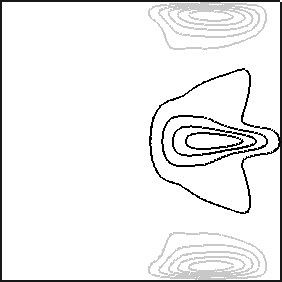}}\hfil
\subfigure[$H_\theta$ at $\theta=0$]{\includegraphics[width=0.3\textwidth]{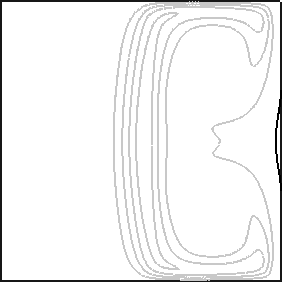}}\hfil
\subfigure[$H_z$ at $\theta=0$]{\includegraphics[width=0.3\textwidth]{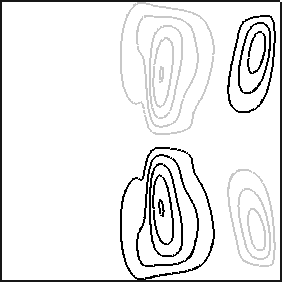}}\\
\subfigure[$H_r$ at $\theta=\frac{\pi}{2}$]{\includegraphics[width=0.3\textwidth]{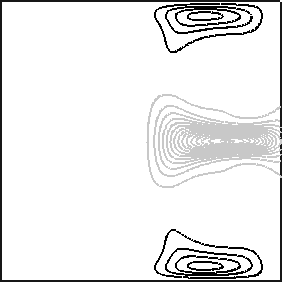}}\hfil
\subfigure[$H_\theta$ at $\theta=\frac{\pi}{2}$]{\includegraphics[width=0.3\textwidth]{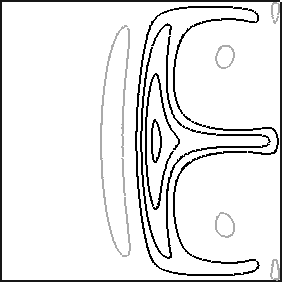}}\hfil
\subfigure[$H_z$ at $\theta=\frac{\pi}{2}$]{\includegraphics[width=0.3\textwidth]{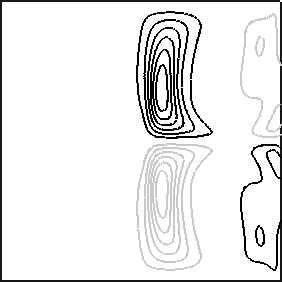}}
\end{minipage}\nopagebreak
\begin{minipage}[c]{0.25\textwidth}%
\includegraphics[width=\textwidth]{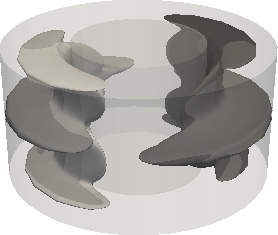} \\
\centerline{(g) isosurface $\|\bHc\|^2$}
\end{minipage} 
\caption{Kinematic dynamo with flow $\calV(\epsilon=8)$ at $\Re=120$, $\Rm=200$.
  Magnetic eigenvector for Fourier mode $m=1$.  Represented in (a) to
  (f) are the radial, azimuthal, and vertical components, normalized
  by the square root of the magnetic energy, in two complementary
  planes, with $0 \le r \le 2, \, -1 \le z \le 1$ (the $z-$axis is on
  the left): for $\theta=0$, $-0.85 \le H_r \le 0.85$ (every 0.17),
  $-0.1 \le H_\theta \le 0.68$ (every 0.17) and $-0.85 \le H-z \le
  0.85$ (every 0.17); for $\theta=\frac{\pi}{2}$, $-1 \le H_r \le
  3.75$ (every 0.25), $-1 \le H_\theta \le 0.15$ (every 0.25) and
  $-1.5 \le H_z \le 1.5$ (every 0.25). Represented in (g) is the
  isosurface $\|\bHc\|^2$ ($7\%$ of maximum value) colored by the
  azimuthal component.  Note the $m=1$ structure.}
\label{fig:H_tcm_eps8}\label{fig:iso_6PCENT_H_tcm_eps8}
\end{figure}
Figure~\ref{fig:H_tcm_eps8} %and~\ref{fig:iso_6PCENT_H_tcm_eps8}
shows the magnetic eigenvector for the Fourier mode $m=1$ at
$\Re=120$, $\Rm=200$ and $\epsilon_{opt}=8$.  This eigenvector is a
rigid wave that rotates in the same direction as the inner cylinder
and top/bottom lids, and its period of rotation is $T \simeq 870$, \ie
more than 66 rotation periods. Upon introducing the equatorial
symmetry operator $\calS_{\text{Z}2} \bH = (H_r,\; H_{\theta},\;
-H_z)(r,\theta,-z)$.  the magnetic field has the following symmetry
property:
\begin{equation}
  \bHc =\calS_{\text{Z}2}\bHc.
\end{equation}
\ie the magnetic field has the same symmetry as the velocity field.
% \begin{figure}[h]
% \includegraphics[width=0.25\textwidth]{\rootfig/iso6PCENT_Emag_eps8_gray_nolegend}
% \caption{Kinematic dynamo at $\Re=120$, $\Rm=200$  and $\epsilon=8$.  
% Represented is the isosurface $|\bH|^2$ ($7\%$ of maximum value)
% colored by the azimuthal component.
% Note the $m=1$ structure characterized
% by one dark blob and one light blob.}
% \label{fig:iso_6PCENT_H_tcm_eps8}
% \end{figure}

\begin{figure}[h]
\includegraphics[width=0.5\textwidth]{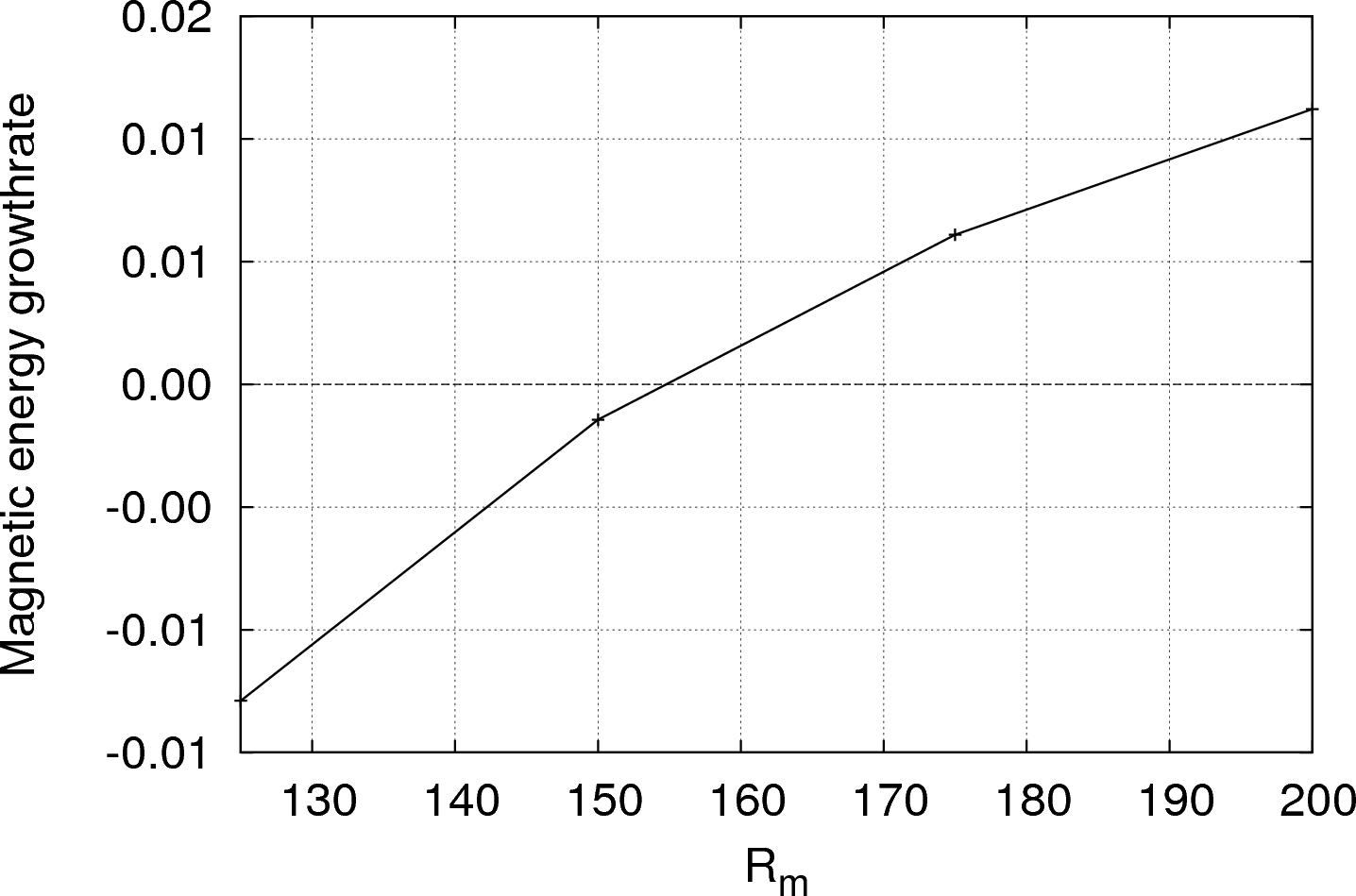}
\caption{Kinematic dynamo of flow $\calV(\epsilon=6.5)$, $\Gamma=2$, $\Re=120$.
  Growthrate of the Fourier mode $m=1$ as a function of $\Rm$.}
\label{fig:sigma_Rm_Re120_EPS6P5}
\end{figure}
We now evaluate the critical magnetic Reynolds number and its minimal
value with respect to $\epsilon$.  We assume that the growthrate
depends smoothly on $\Rm$. Upon inspecting
Figure~\ref{fig:eps_taux_Re120} we see that the growth rate is maximum
for $\epsilon_{opt}\approx 5$ at $\Rm=100$ and for
$\epsilon_{opt}\approx 8$ ate $\Rm=200$. Then, by drawing the line
connecting the two maximum points in Figure~\ref{fig:eps_taux_Re120},
we observe that this line crosses the horizontal line of zero growth
rate in the interval $\epsilon\in [6.5,7]$. We have chosen to explore
the value $\epsilon=6.5$, which gives the poloidal to toroidal ratio
$\Lambda=1.53$.  The growth rate for various magnetic Reynolds
numbers in the range $[125,200]$ has been computed. The results are
shown in Figure~\ref{fig:sigma_Rm_Re120_EPS6P5}.  We have found that
the critical magnetic Reynolds number for the
Fourier mode $m=1$ is $\Rm= 155$ with $\epsilon=6.5$ and $\Lambda=1.53$.

%\begin{figure}[h]
%\epsfig{file=\rootfig/sigma_Rm_Re120_EPS_6P5.eps,width=0.5\textwidth}\hfill
%\caption{Kinematic dynamo. Growthrate of the $m=1$ magnetic energy for Taylor-Couette flow with
%$\Gamma=2$, $\Re=120$ and $\epsilon=6.5$ as a function of $\Rm$.}
%\label{fig:sigma_Rm_Re120_EPS6P5}
%\end{figure}

\subsection{Kinematic dynamo in the  forced Taylor-Couette setup}
\label{sec:kin_forced_TC}
\begin{figure}[h]
\includegraphics[width=0.5\textwidth]{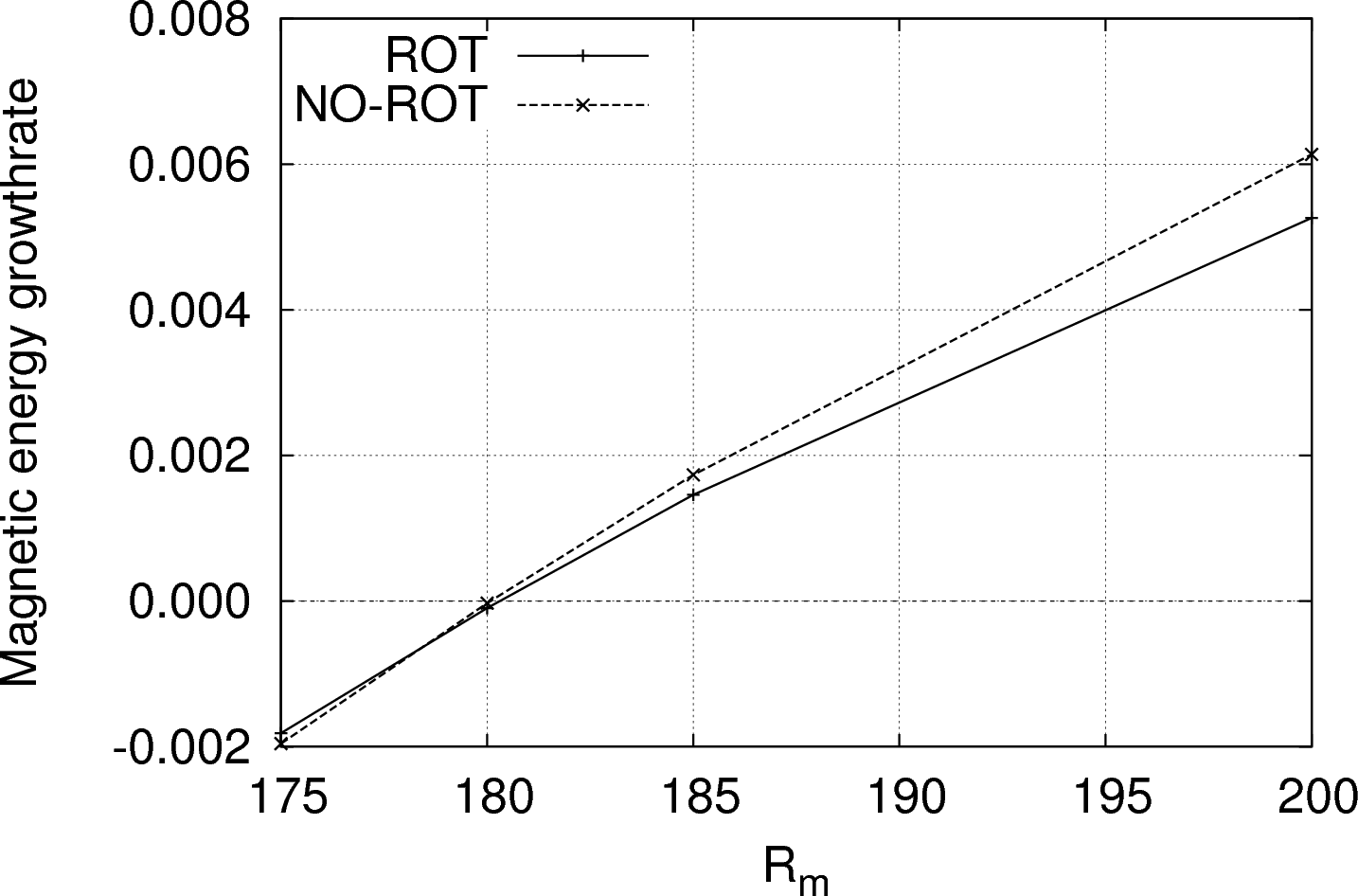}
\caption{Kinematic dynamo with $\calV_I$ flow, $\Gamma=2$ and $\Re=120$
  .  Growthrate of the Fourier mode $m=1$ as a function of $\Rm$. ROT:
  rotating inner core; NO-ROT: non-rotating inner core (but inner wall
  rotates).}
\label{fig:seuil_Re120_Rm}
\end{figure}

We use the steady axisymmetric forced flow, $\calV_{\text{I}}$, at
$\Re=120$ to perform kinematic dynamo computations. The reference
velocity scale is defined in \eqref{def_calU_calV_I} with
$A=2.5$. To determine whether the rotation of
the inner solid core has any impact on the dynamo threshold, we have
compared growth rates when the solid inner core is motionless and when
the inner wall and solid inner core corotate with angular speed
$\Omega_i$. Whether the inner core rotates or not, the magnetic
eigenvector $m=1$ is always the most unstable. We show in
Figure~\ref{fig:seuil_Re120_Rm} the computed growth rates in the two
cases in the range $\Rm\in[175,200]$. The motion of the inner core
does not seem to have a significant influence; in both cases the
dynamo threshold is $\Rm=180$ for the Fourier mode $m=1$.

\begin{figure}[h]
\begin{minipage}{0.70\textwidth}%
\subfigure[$H_r$ at $\theta=0$]{\includegraphics[width=0.31\textwidth]{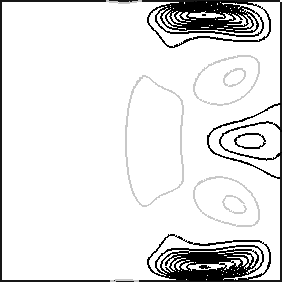}}\hfil
\subfigure[$H_\theta$ at $\theta=0$]{\includegraphics[width=0.31\textwidth]{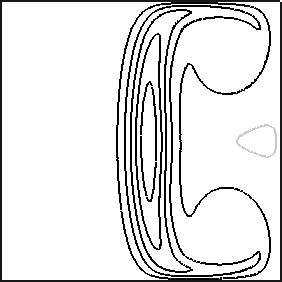}}\hfil
\subfigure[$H_z$ at $\theta=0$]{\includegraphics[width=0.31\textwidth]{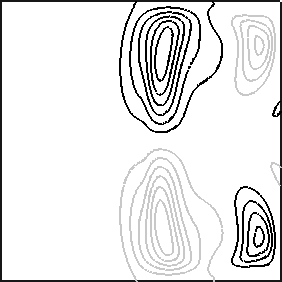}}\\
\subfigure[$H_r$ at $\theta=\frac{\pi}{2}$]{\includegraphics[width=0.31\textwidth]{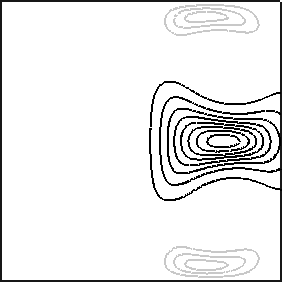}}\hfil
\subfigure[$H_\theta$ at $\theta=\frac{\pi}{2}$]{\includegraphics[width=0.31\textwidth]{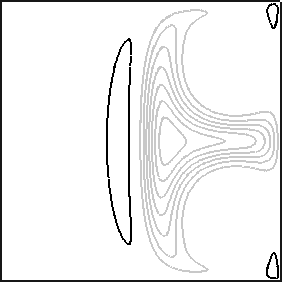}}\hfil
\subfigure[$H_z$ at $\theta=\frac{\pi}{2}$]{\includegraphics[width=0.31\textwidth]{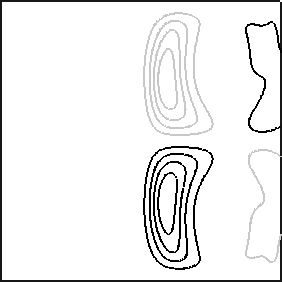}}
\end{minipage}\nopagebreak
\begin{minipage}[c]{0.25\textwidth}%
\subfigure[isosurface $\|\bHc\|^2$]{\includegraphics[width=\textwidth]{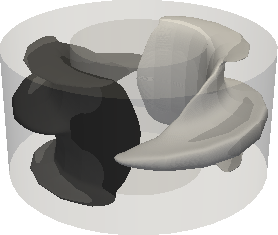}%
\label{fig:iso_14PCENT_H_tcm_Forcage}}%
%\centerline{(g) isosurface $|\bHc|^2$}
\end{minipage} 
\caption{Kinematic dynamo with flow $\calV_I$ at $\Re=120$, $\Rm=200$.
  Magnetic eigenvector for Fourier mode $m=1$.  Represented in (a) to
  (f) are the radial, azimuthal, and vertical components, normalized
  by the square root of the magnetic energy, in two complementary
  planes: for $\theta=0$, $-0.9 \le H_r \le 0.2$ (every 0.1), $-1.4
  \le H_\theta \le 0.35$ (every 0.25) and $-0.6 \le H_z \le 0.6$
  (every 0.1); for $\theta=\pi/2$, $-2.2 \le H_r \le 0.9$ (every
  0.25), $-0.25 \le H_\theta \le 1.75$ (every 0.25) and $-1.25 \le H_z
  \le 1.25$ (every 0.1).  Represented in (g) is the isosurface
  $\|\bHc\|^2$ ($14\%$ of maximum value) colored by the azimuthal
  component.  Note the $m=1$ structure.}
\label{fig:H_tcm_forc_Re120_Rm200}
\end{figure}

The structure of the magnetic eigenvector corresponding to the Fourier
mode $m=1$ with a rotating inner core is shown in
Figure~\ref{fig:H_tcm_forc_Re120_Rm200} (at $\Re=120$, $\Rm=200$). It
is a rigid wave that rotates in the same direction as the inner
cylinder and top/bottom lids with period $T \simeq 120$.
%/workdir/coro/nore/TC_VARGAS/TCM/Re120_Rm200_TCF_AXE_COND/tfin_80/GAYA

\section{Nonlinear dynamo action}
\label{sec:nl_dynamo}
We report in this section on nonlinear dynamo computations done with
the forced Taylor-Couette setup at $\Re=120$ and $\Rm=200$.

\subsection{Description of the setting}
\label{sec:nonlin_forced_TC}
We consider the forced Taylor-Couette setup described in
\S\ref{sec:kin_forced_TC}. The reference
velocity scale is defined in \eqref{def_calU_calV_I} with
$A=2.5$. We perform nonlinear dynamo computations
with the parameters $\Re=120$, $\Rm=200$. The inner core is kept
motionless.  We work with 32 azimuthal modes ($m=0,\ldots,31$), and
the meridional finite element mesh is the same as in the
kinematic runs. The total number of degrees of freedom to be updated
at each time step is $11{,}353{,}104$.  The initial velocity field is
the axisymmetric flow $\calV_{\text{I}}$ that we computed in the Navier-Stokes
regime at $\Re= 120$. The initial magnetic seed is the growing
Fourier mode $m=1$ obtained in the kinematic computations described in
\S\ref{sec:kin_forced_TC} at $\Rm=200 > \Rmc=180$.  

When dynamo action occurs, the magnetic energy grows exponentially
until the Lorentz force is capable of modifying the base flow.  This
transient phase lasts about 5 rotation periods. When the Lorentz force
is strong enough, a new regime settles where the magnetic energy
saturates.  Nonlinear saturation is a slow process that lasts at least
200 rotation periods (see Figure 11 in~\cite{GLLN09}).  Although
SFEMaNS is parallel with respect to the Fourier modes,
%and PETSc~\cite{petsc-web-page} is used to solve in
%parallel the problems in the meridian section, 
the volume of computation required by this type of simulation is such
that we have not been able to explore other kinematic and magnetic
Reynolds numbers within the resources allocated to this project. The
nonlinear run presented in this section used about 15600 cumulated CPU
hours with 32 processors on an IBM Power 6 cluster.

\subsection{Time evolution of the energy}
\begin{figure}[h]
\subfigure[Kinetic (-0.37) and magnetic energies]{\includegraphics[width=0.5\textwidth]{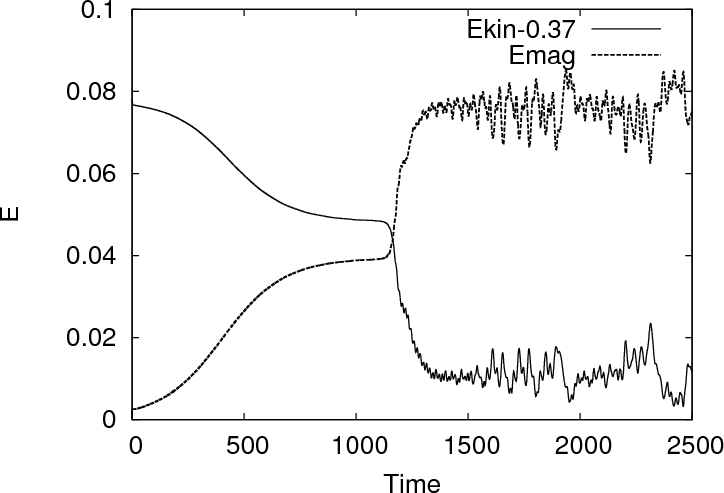}}\hfil
%\subfigure[Kinetic energy]{\includegraphics[width=0.33\textwidth]{\rootfig/Ecin_fini}}\hfil
%\subfigure[Magnetic energy]{\includegraphics[width=0.33\textwidth]{\rootfig/Emag_fini}}\hfil
\subfigure[Splitting of magnetic energy]{\includegraphics[width=0.5\textwidth]{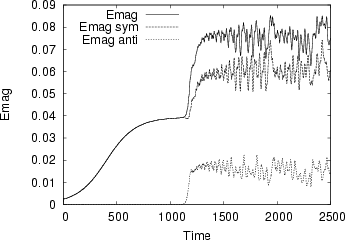}}
\caption{Nonlinear dynamo in the forced Taylor-Couette setup.  (a)
  Time evolution of kinetic (-0.37) and magnetic energies in the conducting
  region $0 \le r \le R_o$ and $-\Gamma/2 \le z \le \Gamma/2 $. Panel
  (b) shows the symmetric and anti-symmetric components of the
  magnetic energy.}
\label{fig:en_time_fini}
\end{figure}

The time evolutions of the kinetic and magnetic energies are reported
in Figure~\ref{fig:en_time_fini}(a-b), where the kinetic and magnetic
energies are defined as follows: $\frac12 \int_{\Omegacf} \|\bu\|^2
\dif \bx$, $\frac12 \int_{\Omegac} \|\bHc\|^2 \dif \bx$, respectively.
From $t=0$ to $t=500$ (first transition), the kinetic energy decreases
and the magnetic energy grows exponentially with a growthrate similar
to that of the kinematic dynamo.  Then both the magnetic and the
kinetic energies seem to saturate in a first nonlinear regime,
$500\lesssim t \lesssim 1100$. During the first transition, the fluid
flow loses the axial symmetry, $m=0$, and the magnetic field loses the
symmetry associated with the Fourier mode $m=1$.  The flow being
forced by the Lorentz force $(\ROT \bHc){\times}\bHc$, the velocity
thereby acquires a contribution on the Fourier mode $m=2$.  The
magnetic field being deformed by the action of the induction term
$\bu{\times}\bHc$ acquires a contribution on the Fourier mode
$m=3$. The cascade of nonlinear couplings generate even velocity
modes, $m=0,2,\ldots$, and odd magnetic modes, $m=1,3,\ldots$ During
this transitional phase that consists of populating the Fourier modes,
the axisymmetry of the velocity field is broken but the equatorial
(mid-plane) symmetry is preserved for both the velocity and the
magnetic field.

%/workdir/coro/nore/TC_VARGAS/TCMHD/Re120_Rm200_TCF_AXE_COND/M32

\begin{figure}[h]
%\subfigure[Magnetic energy, $m=0,1,2,3$]{\includegraphics[width=0.32\textwidth]{\rootfig/Emag_modes0123}}\hfil
%\subfigure[Kin. $m= 1,\,3$; Mag. $m=0,\,2$]{\includegraphics[width=0.32\textwidth]{\rootfig/Ekin_mag_m_small}}\hfil
%\subfigure[Magnetic energy, $m=0,1,30,31$]{\includegraphics[width=0.32\textwidth]{\rootfig/Emag_modes013031}}
\subfigure[Magnetic energy, $m=0,1,2,3$]{\includegraphics[width=0.33\textwidth]{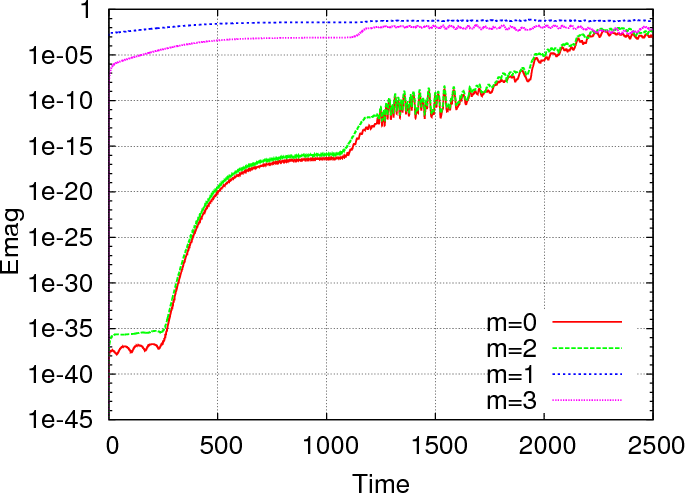}}\hfil
\subfigure[Kin. $m= 1,\,3$; Mag. $m=0,\,2$]{\includegraphics[width=0.33\textwidth]{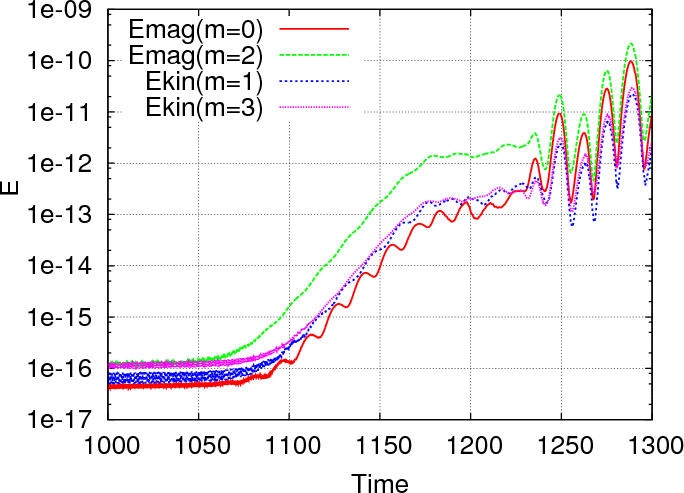}}\hfil
\subfigure[Magnetic energy, $m=0,1,30,31$]{\includegraphics[width=0.33\textwidth]{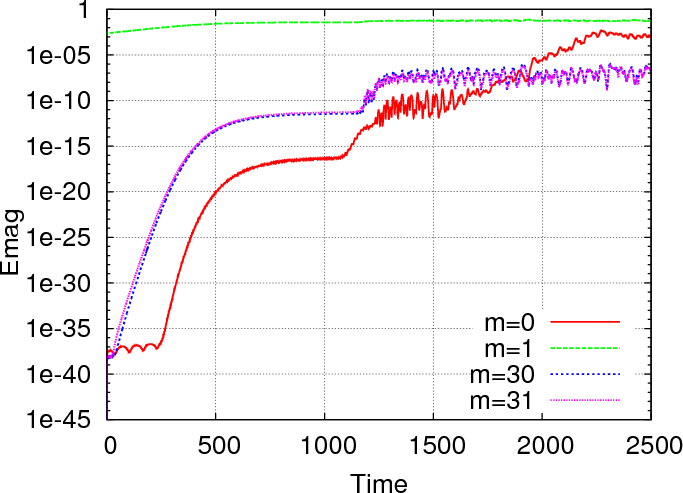}}%
\caption{Time evolution of different modal energies; (a) magnetic
  energies $m=0, \,2$ (bottom curves) and $m= 1, \,3$ (top curves);
  (b) kinetic energies $m= 1, \,3$ and magnetic energies $m=0, \,2$;
  (c) magnetic energies $m=1$ (top curve), $m=30, \,31$ (middle
  curves) and $m=0$.}
  \label{fig:Ekin_mag_modes}
\end{figure}

A second nonlinear transition starts at $t=1100$ and lasts until $t=
1175$. In this time interval the magnetic energy increases and the
kinetic energy decreases. This change of behavior is due to the
breaking of the equatorial symmetry. This phenomenon is well
illustrated by computing the energy of the symmetric part,
$\frac12(\bHc + \calS_{\text{Z}2}\bHc)$, and anti-symmetric part,
$\frac12(\bHc - \calS_{\text{Z}2}\bHc)$, of the magnetic field.  The
time evolution of these two quantities is shown in
Figure~\ref{fig:en_time_fini}(c).  The equatorial symmetry breaking is
driven by the small even azimuthal modes of the magnetic field as can
be seen on Figure~\ref{fig:Ekin_mag_modes} (a)(b), especially the
magnetic mode $m=2$.

In the time interval $1175 \le t \le 1600$, the system enters a
second nonlinear regime characterized by large fluctuations and a
dynamics dominated by the large Fourier modes.  Between $t=1600$ and
$t= 2200$, we observe a third transition during which the small even
modes of the magnetic field increase again until they reach the final
saturated state. A third and final nonlinear regime settles beyond
$t=2200$. Figure~\ref{fig:Ekin_mag_modes}(c) shows that the large
Fourier modes, exemplified by $m=30$ and $m=31$, basically fluctuate
within some asymptotic range for $t> 1250$, whereas the small even
modes grow until they become energetically significant.

\begin{figure}[h]
\subfigure[$H_\theta(1.2,0,-0.5,t)$]{\includegraphics[width=0.4\textwidth]{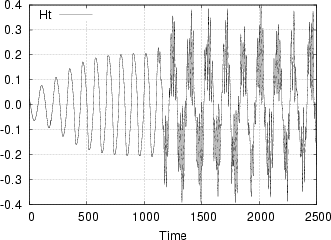}\label{fig:anemo_Ht_pt4_a}}\hfil
\subfigure[$H_\theta(1.2,0,-0.5,t)-H_\theta(1.2,0,0.5,t)$]{\includegraphics[width=0.4\textwidth]{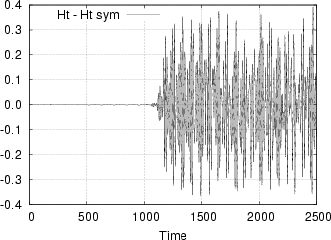}\label{fig:anemo_Ht_pt4_b}}
\caption{Time series of (a) $H_\theta$ at $(r=1.2, \theta=0, z=-0.5)$
  and (b) twice the anti-symmetric part of $H_\theta$ at $(r=1.2,
  \theta=0, z=-0.5)$.}
\label{fig:anemo_Ht_pt4}
\end{figure}

We show in Figure~\ref{fig:anemo_Ht_pt4_a} the time series of the
azimuthal component of the magnetic field at the point
$(r=1.2,\theta=0,z=-0.5)$. The envelop of the signal first grows
then reaches a maximum range. The period of the signal is $T \approx
112.5$; this corresponds to a wave that rotates in the same direction
as the inner cylinder and top/bottom lids. More frequencies appear
beyond $t=1100$; the signal is the superposition of an oscillation of
period $T \approx 150$ and a modulation of period $T_{\text{mod}} \approx
17$, which happens to be of the same order as the wall rotation period
$T_{\text{lids}}= 2 \pi/0.55 \approx 11.4$. The breaking of the equatorial
symmetry is measured by monitoring the anti-symmetric part of the
magnetic field. We show in Figure~\ref{fig:anemo_Ht_pt4_b} the time
evolution of twice the anti-symmetric part of $H_\theta$ at
$(r=1.2,\theta=0,z=-0.5)$.

\subsection{Spatial structure of the dynamo}
\begin{figure}[h]
\begin{minipage}[b]{0.24\textwidth}Isosurface $\|\bHc\|^2$ ($25\%$ of maximum value)\end{minipage}%
\subfigure[ $t=1000$]{\includegraphics[width=0.24\textwidth]{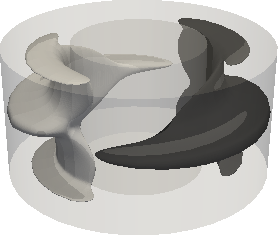}}\hfil
\subfigure[ $t=1400$]{\includegraphics[width=0.24\textwidth]{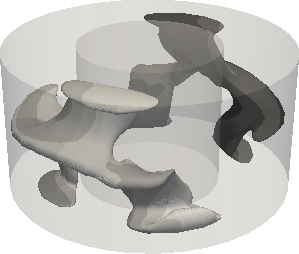}}\hfil
\subfigure[ $t=2500$]{\includegraphics[width=0.24\textwidth]{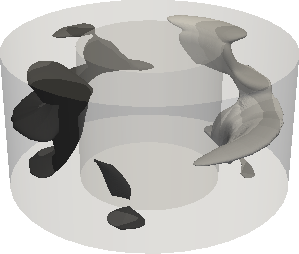}}
\begin{minipage}[b]{0.24\textwidth}Isosurface $\|\bV\|^2$ ($25\%$ of maximum value)\end{minipage}%
\subfigure[ $t=1000$]{\includegraphics[width=0.24\textwidth]{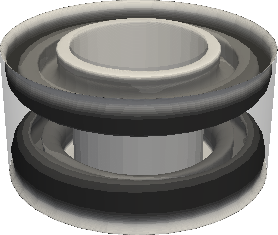}}\hfil
\subfigure[ $t=1400$]{\includegraphics[width=0.24\textwidth]{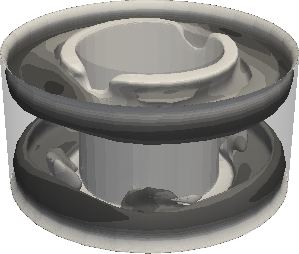}}\hfil
\subfigure[ $t=2500$]{\includegraphics[width=0.24\textwidth]{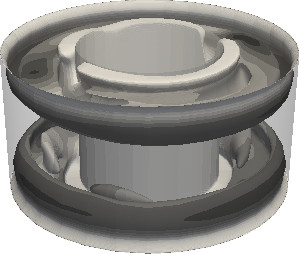}}
\begin{minipage}[b]{0.24\textwidth}Isosurface $\|\bV\|^2$ without the
  axisymmetric mode $m\,{=}\,0$ ($10\%$ of maximum value)\end{minipage}%
\subfigure[$t=1000$]{\includegraphics[width=0.24\textwidth]{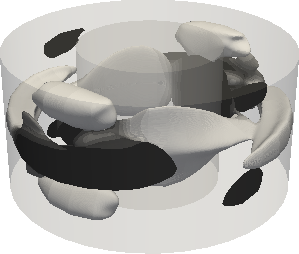}}\hfil
\subfigure[$t=1400$]{\includegraphics[width=0.24\textwidth]{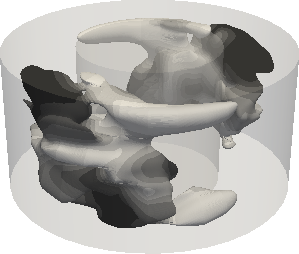}}\hfil
\subfigure[$t=2500$]{\includegraphics[width=0.24\textwidth]{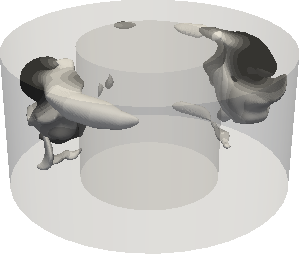}}
\caption{Nonlinear dynamo at $t=1000$, $t=1400$ and $t=2500$, for
  $\Re=120$, $\Rm=200$: (a-c) isosurface of $\|\bHc\|^2$ ($25\%$ of
  maximum value); (d-f) isosurface of $\|\bV\|^2$ ($25\%$ of maximum
  value); (g-i) isosurface of $\|\bV\|^2$ without the axisymmetric
  mode ($10\%$ of maximum value); color scale proportional to azimuthal component.}
\label{fig:iso_TCf_Re120_Rm200}
\end{figure}

To have a better understanding of the structure of the velocity and
magnetic fields during the three nonlinear regimes identified above,
we show in Figure~\ref{fig:iso_TCf_Re120_Rm200} the isosurfaces of the
magnetic and kinetic energies at $t=1000$, $t=1400$, and
$t=2500$. Since the axisymmetric velocity mode $m=0$ is dominant (see
Figure~\ref{fig:iso_TCf_Re120_Rm200}(d-f)), we show in
Figure~\ref{fig:iso_TCf_Re120_Rm200}(g-i) the kinetic energy without
its axisymmetric contribution to better distinguish the fine
structures.  The magnetic field at $t=1000$ is dominated by the odd
azimuthal modes (see Figure~\ref{fig:iso_TCf_Re120_Rm200}(a)) and
resembles the eigenvector shown in
Figure~\ref{fig:iso_14PCENT_H_tcm_Forcage}. The velocity field is
composed of even modes as can be seen on
Figure~\ref{fig:iso_TCf_Re120_Rm200}(g), see \eg the two dark
structures that are diametrically opposed. Similar spatial
distributions are observed at time $t=1400$ with the addition of
smaller scales.  At time $t=2500$ the small scale modes are even more
apparent, and the magnetic and velocity structures in the fluid domain
are more deformed.

\begin{figure}[h]
\subfigure[$0 \le t \le 2500$]{\includegraphics[width=0.45\textwidth]{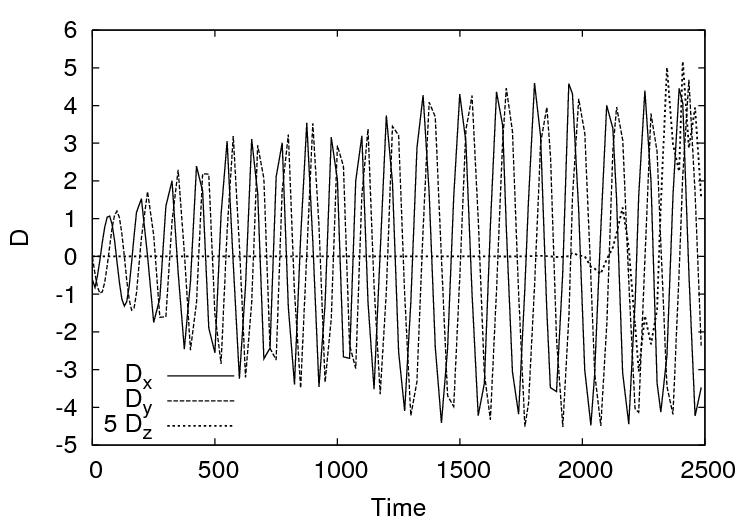}}\hfil
\subfigure[Zoom in $2410 \le t \le 2570$]{\includegraphics[width=0.45\textwidth]{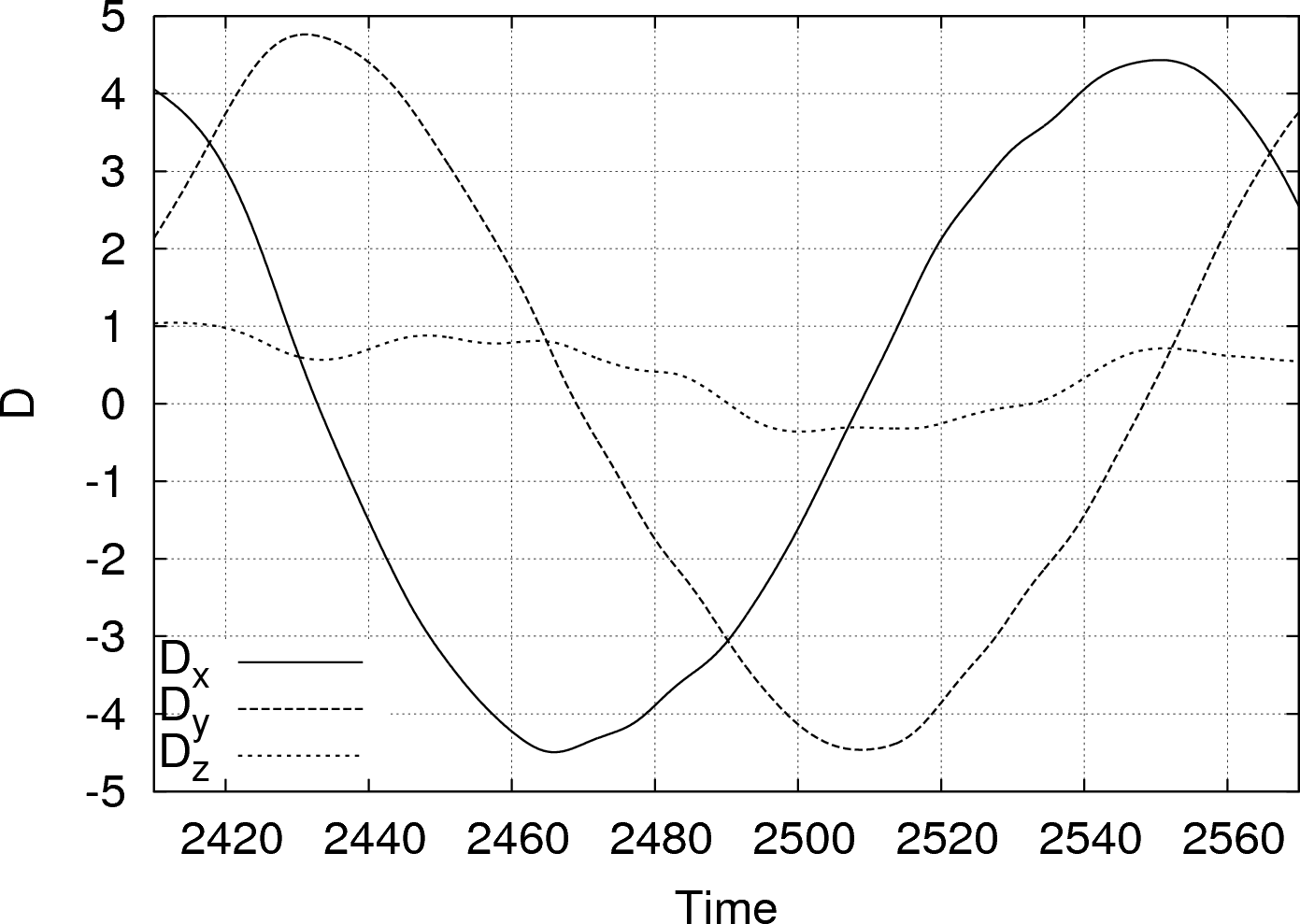}}
\caption{Time evolution of the dipolar magnetic moment. Note that the vertical component has been
multiplied by 5 in the left panel.}
  \label{fig:dip}
\end{figure}
\begin{figure}[h]
\subfigure[From above]{\includegraphics[width=0.45\textwidth]{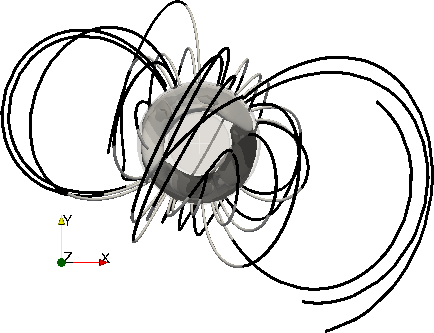}}\hfil
\subfigure[From a perspective point of view]{\includegraphics[width=0.45\textwidth]{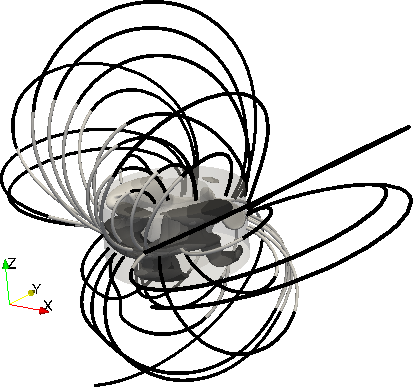}}
\caption{Isosurface of $\|\bHc\|^2$ ($25\%$ of maximum value) and
  magnetic field lines at $t=2500$.}
  \label{fig:Re120_Rm200_t2500}
\end{figure}

To characterize the long distance influence of the magnetic field, we
have recorded the time evolution of the magnetic dipole defined by
$\bD= \int_{\Omegac} \br {\times} (\ROT\bHc) \dif \bx$.
Figure~\ref{fig:dip} shows the time series of the three Cartesian
components of the magnetic dipole in the time interval $0\le t \le
2500$. During the first two transitions and nonlinear regimes, \ie
$0\le t\le 1600$, the dipolar moment is purely equatorial and rotates
at the same frequency as the magnetic field.  The axial moment starts
to grow at the beginning of the third transition ($t > 1600$) and
changes sign several times afterward (note that the time series in (a)
is under sampled in this range). Figure~\ref{fig:dip}(b) presents a
zoom of the time evolution of $D_z$ showing two reversals.  The
magnitude of the quadrupolar moments (data not shown) stay below
$10^{-3}$ until $t=1175$, then increase and saturate to values four
times smaller than the magnitude of the dipolar moment. Magnetic field
lines in the vacuum at $t=2500$ show a pattern characteristics of an
equatorial dipole (see Figure~\ref{fig:Re120_Rm200_t2500}).

\section{Concluding remarks}  \label{Sec:remarks}
We have numerically demonstrated that pure viscous driving by
smooth rotating walls in a short Taylor-Couette setup does not lead to
dynamo action for $\Rm \le 200$ since the poloidal to toroidal ratio of the velocity
field is too small. An adjustment of the poloidal to toroidal ratio is
needed to achieve dynamo action in the kinematic regime. We have
implemented an ad hoc body force to produce a poloidal to
toroidal ratio that is of the same order as what is needed in the
kinematic simulations to trigger the dynamo action. This force
may be thought of as a model for the action of blades fixed to the
static or to the rotating walls/lids at convenient angles.
%and with an adapted width to be optimized in a water model.
This force has also the same symmetry properties as the
geodynamo, \ie the $\text{SO}(2){\times}\text{Z}2$ symmetry
(axisymmetry and equatorial symmetry).  The critical magnetic Reynolds
number of this setup based on the inner cylinder speed $\mu_0 \sigma
\Omega_i R_i^2=0.55{\times}180=99$ is in the range of what has been
obtained in the kinematic studies of~\cite{Dudley89} in a spherical
container with the same symmetries. This magnetic Reynolds number is
also comparable to what has been reported in~\cite{GLLN09} for
Taylor-Couette simulations in vessels of larger aspect ratios and with
pure viscous driving.

A nonlinear simulation has been performed at $\Re=120$, $\Rm=200$ over
225 rotation periods. In the early linear phase of the dynamo, the
external field is dominated by an equatorial rotating dipole.  In the
established nonlinear regime, an axial axisymmetric component is
excited and exhibits reversals.  The relation between the main flow
parameters of the time-dependent angle formed by the dipole and the
rotation axis calls for further investigations, since it is a basic
feature of observed planetary dynamos.

Exploring the feasibility of an experimental fluid dynamo based on the
present design will require expertise from many different experimental
and numerical fields~\cite{LN08}.  To achieve a magnetic Reynolds
number equal to $100$ in a flow of liquid sodium requires that the
kinematic Reynolds number be of order $10^7$.  It is well known that
such a value corresponds to a highly turbulent flow that can be
studied only in experimental facilities, since it is far beyond the
capacity of direct numerical simulations. The objective of such
experiments should be to recover optimized poloidal and toroidal
components after time averaging, which presumably would guide the
design of the blades fixed to the endwalls. These experiments would
also inform about the power requirements. Using a standard rotation
frequency of $50$~Hz, a magnetic Reynolds number of $100$ can be
obtained in liquid sodium at $150^{\text{o}}$C with an inner radius of
approximatively $18$~cm and an outer radius and height of
$36$~cm. This seems feasible since these dimensions are not far from
those of the Cadarache experiment~\cite{Monchaux07}. We conjecture
however that the power required by this experiment at a given rotation
frequency should be smaller, since the turbulence rate induced by
co-rotating lids/impellers should be smaller than that of
counter-rotating lids/impellers.  A dynamo facility presenting
similarities with the present proposal is currently investigated by
Colgate and collaborators~\cite{Colgate2011}. Their MHD device uses
also a Taylor-Couette forcing in a short cylindrical container with
size and targeted magnetic Reynolds number similar to those studied in
the present paper. There are however differences: the flow in their
experiment forms an outwards jet in the equatorial plane and is driven
by viscous stresses only.  More detailed comparisons of the respective
merits of both designs should certainly be instructive.

%(Rm=muo sigma0 2 pi N R^2 =10 2 pi *50 *.2^2)

\section*{Acknowledgments}
JLG is thankful to University Paris Sud 11 for constant support over
the years.  The computations were carried out on the IBM Power 6
cluster of Institut du D\'eveloppement et des Ressources en
Informatique Scientifique (IDRIS) (project \# 0254).  \vspace{2cm}
\bibliographystyle{plain} \bibliography{biblio}
%\bibliography{/u/guermond/GUERMOND/BIBLIO/biblio}

\end{document}